%% file: axioms_safem.tex
\documentclass{siamltex1213}

\usepackage[T1]{fontenc}
\usepackage{amsmath}     %eqref
\usepackage{amssymb}     
\usepackage{esint} %fuer \fint: \usepackage{txfonts} oder {pxfonts}
\usepackage{pgf}
\usepackage{tikz}
\usetikzlibrary{decorations.pathmorphing}
\usetikzlibrary{decorations.pathreplacing}
\usetikzlibrary{positioning}
\usetikzlibrary{shapes}
\usetikzlibrary{arrows}
\usetikzlibrary{patterns}
\usetikzlibrary{fadings}
\usetikzlibrary{plotmarks}
\usetikzlibrary{calc}
\usetikzlibrary{intersections}
\tikzstyle{every picture}+=[font=\footnotesize]
\usepackage{paralist}

\usepackage{bbm}
\usepackage{latexsym}           %benoetigt fuer \Box am Beweisende
\usepackage{subfigure}			%fuers einfuegen von eps Grafiken
\usepackage{enumerate}
\usepackage{enumitem}

\setlist{noitemsep, topsep=0.8ex, partopsep=0pt%, parsep=0pt, itemsep=0pt
	, leftmargin=3em}
\setlist[1]{labelindent=\parindent}

\newlist{axioms}{enumerate}{1}
\setlist[axioms]{font=\bfseries}

\newlist{alphenum}{enumerate}{1}
\setlist[alphenum]{label=\textbf{(\alph*)}, leftmargin=4em}

\newlist{alphienum}{enumerate}{1}
\setlist[alphienum]{label=\textit{(\alph*)}}

\newlist{romanenum}{enumerate}{1}
\setlist[romanenum]{label=\textit{(\roman*)}}

\newlist{romaninenum}{enumerate*}{1}
\setlist[romaninenum]{label=\textit{(\roman*)}}
\usepackage[vlined]{algorithm2e}
\SetKwIF{If}{ElseIf}{Else}{if}{}{else if}{else}{endif}
\SetKwFor{For}{for}{}{endfor}
\usepackage[noabbrev, capitalise]{cleveref}
\usepackage{tabu}
\crefname{equation}{\unskip}{\unskip}
\creflabelformat{equation}{#2(#1)#3}
        \SetFuncSty{textsc}
        \SetKwFunction{frun}{Run}
        \SetKwHangingKw{arun}{\frun}
        \SetKwFunction{fset}{Set}
        \SetKwHangingKw{aset}{\fset}
        \SetKwFunction{fselect}{Select}
        \SetKwHangingKw{aselect}{\fselect}
        \SetKwFunction{fcompute}{Compute}
        \SetKwHangingKw{acompute}{\fcompute}
        \SetKwFunction{fsolve}{Solve}
        \SetKwHangingKw{asolve}{\fsolve}
        \SetKwFunction{festimate}{Estimate}
        \SetKwHangingKw{aestimate}{\festimate}
        \SetKwFunction{fmark}{Mark}
        \SetKwHangingKw{amark}{\fmark}
        \SetKwFunction{frefine}{Refine}
        \SetKwHangingKw{arefine}{\frefine}
        \SetKwFunction{compute}{Compute}
        \SetKwFunction{set}{Set}
\newtheorem{remark}{Remark}
\input{command}

\allowdisplaybreaks[0] %ermoeglicht seitenumbrueche in Formeln 
\newenvironment{algof}%
{\algorithm}%
   {\endalgorithm}% Figure/float layout
   
\title{Axioms of adaptivity for separate marking}
\author{C.~Carstensen\footnotemark[1] \footnotemark[2] \and H.~Rabus\footnotemark[1]}
\begin{document}
\maketitle
\renewcommand{\thefootnote}{\fnsymbol{footnote}}
\footnotetext[1]{Department of Mathematics, Humboldt-Universit\"at zu Berlin, Unter den Linden 6,
        10099 Berlin, Germany.   Email cc@math.hu-berlin.de and rabus@math.hu-berlin.de
}
\footnotetext[2]{The work of the first author is partly supported by DFG SPP 1749 \textit{Reliable Simulation Techniques in Solid Mechanics -- Development of Non-standard Discretization Methods, Mechanical and Mathematical Analysis}}
	\renewcommand{\thefootnote}{\arabic{footnote}}

\pagestyle{myheadings}
\thispagestyle{plain}
\markboth{C.~Carstensen and H.~Rabus}{Axioms of adaptivity for separate marking}
%\subjclass{Primary 65N12, 65N15, 65N30, 65N50, 65Y20}

%%%%
% abstract
%%%%
\begin{abstract}
Mixed finite element methods with flux errors in $H(\ddiv)$-norms and div-least-squares 
finite element methods require a separate marking strategy in obligatory adaptive mesh-refining.
The refinement indicator $\sigma^2(\T,K)=\eta^2(\T,K)+\mu^2(K)$ of a finite 
element domain $K$ in an admissible triangulation $\T$ consists of some 
residual-based error estimator $\eta(\T,K)$ with some reduction property under local 
mesh-refining and some data approximation error $\mu(K)$. Separate marking means either D\"orfler
marking if $\mu^2(\T) \leq \kappa \eta^2(\T)$ or otherwise an optimal data approximation algorithm runs with
controlled accuracy as established in \cite{CR09, safem2015}. 

The axioms are abstract and sufficient conditions on the estimators $\eta(\T,K)$ and data
approximation errors $\mu(K)$ for optimal asymptotic convergence rates. The enfolded set of axioms simplifies \cite{CFP14} for collective marking,
treats separate marking established for the first time in an abstract framework, generalizes \cite{CCP-lsfem}
for least-squares schemes, and extends \cite{CR09} to the mixed FEM with flux error control in $H(\ddiv)$.
\end{abstract}

\input{intro}
\input{prelim}

\input{remarks}
\input{optimality}

\input{applications}
\bibliographystyle{alpha}
\bibliography{biblio}

\end{document}

%% file: command.tex
\DeclareMathOperator{\D}{D}

\DeclareMathOperator{\curl}{curl}

\DeclareMathOperator{\ddiv}{div}

\DeclareMathOperator{\Hdiv}{\mathrm H( div, \Omega )}

\DeclareMathOperator{\osc}{osc}

\DeclareMathOperator{\Tol}{Tol}

\newcommand{\T}{\mathcal{T}}
\renewcommand{\P}{\mathcal{P}}
\newcommand{\hT}{\hat{\mathcal{T}}}

\newcommand{\hTl}[1]{\hat{\mathcal{T}}_{#1}}
\newcommand{\tTl}[1]{\tilde{\mathcal{T}}_{#1}}
\newcommand{\TT}[1]{\mathbb T\left(#1\right)}
\newcommand{\PP}[1]{\mathbb P\left(#1\right)}
\newcommand{\TO}{\mathbb T}
\newcommand{\POO}{\mathbb P}

\DeclareMathOperator{\E}{\mathcal{E}}
\newcommand{\M}{\mathcal{M}}
\newcommand{\Tl}[1]{\mathcal T_{#1}}

\newcommand{\Ml}[1]{\mathcal M_{#1}}

\newcommand{\NO}[0]{\mathbb{N}_0} % natural numbers

\newcommand{\norm}[2]{\left\lVert #1 \right\rVert_{#2}}
\newcommand{\restrict}[2]{\left. #1 \right\lvert_{#2}}
\newcommand{\abs}[1]{\left\lvert #1 \right\rvert}

\newcommand{\LO}[0]{L^2(\Omega)}

\newcommand{\Lstab}[0]{\ensuremath{\Lambda_{\mathrm{1}}}}
\newcommand{\Lred}[0]{\ensuremath{\Lambda_{\mathrm{2}}}}
\newcommand{\Lqo}[0]{\ensuremath{\Lambda_{\mathrm{4}}}}

\newcommand{\LqoE}[0]{\ensuremath{\Lambda_{\mathrm{4}(\varepsilon)}}}
\newcommand{\LqoEs}[0]{\ensuremath{\Lambda_{\mathrm{4}(\varepsilon')}}}
\newcommand{\LdRel}[0]{\ensuremath{\Lambda_{\mathrm{3}}}}
\newcommand{\LhdRel}[0]{\ensuremath{\widehat{\Lambda}_{\mathrm{3}}}}
\newcommand{\whM}[0]{\ensuremath{\widehat{M}}}
\newcommand{\LoptData}[0]{\ensuremath{\Lambda_{\mathrm{5}}}}
\newcommand{\LsubAdd}[0]{\ensuremath{\Lambda_{\mathrm{6}}}}
\newcommand{\Lqmeta}[0]{\ensuremath{\Lambda_{\mathrm{8}}}}
\newcommand{\LoptA}[0]{\ensuremath{\Lambda_{\mathrm{9}}}}

\newcommand{\Lqm}[0]{\ensuremath{\Lambda_{\mathrm{7}}}}
\newcommand{\Lref}[0]{\ensuremath{\Lambda_{\mathrm{ref}}}}
\newcommand{\Lopt}[0]{\ensuremath{\Lambda_{\mathrm{opt}}}}
\newcommand{\Lcl}[0]{\ensuremath{\Lambda_{\mathrm{BDdV}}}}
\newcommand{\LII}[0]{\ensuremath{\Lambda_{\mathrm{12}}}}
\newcommand{\La}[0]{\ensuremath{\Lambda}}

\newcommand{\LestRed}[0]{\ensuremath{R_\ell}}
\newcommand{\estJ}[1]{\eta_{#1}}
\newcommand{\estV}[1]{\mu_{#1}}

\newcommand{\volO}[1]{\norm{h_{#1}f}{\LO}}

\newcommand{\cafem}{\textsc{cafem}\xspace}
\newcommand{\safem}{\textsc{safem}\xspace}
\newcommand{\safems}{\textsc{safem}s\xspace}
\newcommand{\afem}{\textsc{afem}\xspace}
\newcommand{\afems}{\textsc{afem}s\xspace}

\newcommand{\fem}{\textsc{fem}\xspace}

\def\pRT{p_{RT}}
\def\qRT{q_{RT}}
\def\hpRT{\widehat{p_{RT}}}
\def\huRT{\widehat{u_{RT}}}
\def\hqRT{\widehat{q_{RT}}}

\def\midhT{\operatorname{mid}(\hT)}
\def\huCR{\widehat{u_{CR}}}
\def\hvCR{\widehat{v_{CR}}}

\newcommand{\Refine}{{\scshape Refine}\xspace}

\newcommand{\Approx}{{\scshape Approx}\xspace}

%% file: intro.tex
\section{Introduction}
The convergence analysis of adaptive finite element methods (\afems) with collective marking for some total error estimator (called \cafem below) is reformulated in an abstract setting in \cite{CFP14}. Therein four axioms describe elementary properties of the total error estimator that are sufficient for optimal convergence rates. 
Standard adaptive schemes are based on a total error estimator and collective marking on each level outlined in pseudo code as follows. 
\medskip

\begin{minipage}{0.9\textwidth}
	\textbf{CAFEM$(\theta, \Tl 0)$}
	\\
	\input{algCAfem}

\end{minipage}

This paper simplifies the axioms from \cite{CFP14},  also works without the concept of nonlinear approximation classes \cite{BDD04,Stev07,CKNS07} and so avoids  any notion of efficiency.
The recent comprehensive  a~posteriori error analysis in \cite{ccdpas2015} provides  
an efficient and reliable  control in natural norms: the error in the flux in $H(\ddiv,\Omega)$ 
and the error in the displacements in $L^2(\Omega)$. 
The focus of this paper is on separate marking (\safems), a modification of the standard \afem: D\"orfler marking is applied if the estimated error dominates the data approximation error, while an optimal data approximation is performed otherwise --- outlined in pseudo code as follows.
\medskip

\begin{minipage}{0.9\textwidth}
	\textbf{SAFEM$(\theta_A, \kappa, \rho_B, \Tl 0)$}
	\\
	\input{algSAfem}

\end{minipage}

The algorithm \safem combines ideas from \cite{mfemBeckerMao08,CR09,safem2015} and distinguishes two Cases (A) and (B), where the refinement is with respect to the dominant refinement indication $\eta_\ell^2$ or $\mu_\ell^2$. The refinement in Case (B) depends on the data approximation error and is independent of the discrete solution. This allows for any optimal algorithm for data approximation with respect to the error functional $\mu^2:K \rightarrow \mathbb R$ for $K\subseteq \Omega\subseteq \mathbb R^n$, i.e.\ the output $\Tl{\Tol}=\texttt{appx}(\Tol, \mu(K):K\in\Tl 0)$ is expected to satisfy
\begin{align*}
	\mu^2(\Tl \Tol)&\leq \Tol,\\
	\abs{\Tl \Tol} - \abs{\Tl 0} &\leq \LoptData \Tol^{-1/s}.
\end{align*}
The analysis for \afems based on collective marking as in \cite{CFP14} is included when $\sigma^2(\T,\bullet)=\eta^2(\T,\bullet)+\mu^2(\T,\bullet)$ replaces $\eta^2(\T,\bullet)$ in Case (A) and the refinement indicator in Case (B) vanishes.

Optimal convergence rates for the estimators follow from  axioms (A1)-(A4) 
generalized from \cite{CFP14} and (B1)-(B2) for optimal  data approximation with quasimonotonicity (QM). 
The subroutine  \texttt{appx}  in \safem can be realized by some D\"orfler marking (similar to the algorithm in \cite{mfemBeckerMao08}) or by the algorithm \Approx  from \cite{BDD04,BD04}  (applied in \cite{CR09,safem2015}). The flexibility in the data reduction allows
applications of \safem to  problems with data approximation terms that do \textit{not} satisfy an
estimator reduction property but quasimonotonicity.  Two model examples illustrate this in the present paper: mixed \fem with flux error estimation in $H(\ddiv)$ rather then $L^2(\Omega)$ \cite{CR09} and a least-squares \fem problem from \cite{CCP-lsfem}. Further applications of the present version of the axioms on \safem shall
appear in the near future \cite{lsfemBC,BCStarke}.

The remaining parts of this paper are organised as follows. Section \ref{sec:axioms} presents more details on \safem and guides the reader through the conditions in (A1)-(A4) and (B1)-(B2) for the refinement indicators $\eta$ and $\mu$ and asserts the optimal convergence rate of \safem in \cref{thm:safem}. 
A collection of remarks follows in Section \ref{sec:remarks} before Section \ref{sec:optimality} presents the proofs. Sections \ref{s:applmfem}-\ref{s:appllsfem} contain the verification of the axioms for two examples, where separate marking is obligatory for optimal adaptive mesh-refinement.
The main novel contribution in Section \ref{s:applmfem} is the proof of a discrete version  (A3) of \cite{ccdpas2015}.

The notation $A \lesssim B$ abbreviates $A \leq CB$ for some positive generic constant $C$, which depends only on the initial triangulation $\Tl 0$ and on the universal constants in the axioms; while $A\approx B$ abbreviates $A\lesssim B \lesssim A$.
Throughout this paper standard notation of Lebesgue and Sobolev spaces and their norms applies. The modulus sign 
$|\bullet|$ denotes the Euclidean length as well as the counting measure, e.g., $|\mathcal{M}|$ is the cardinality of 
$\mathcal{M}$ and equals the number of elements
in a triangulation $\mathcal{M}$ (or a subset thereof). 

%% file: algCAfem.tex
\begin{algof}
	\DontPrintSemicolon
	\For{$\ell=0,1,\ldots$}{
		\acompute{ $\sigma_\ell(K)$ for all $K\in \Tl \ell$
		}
		{$\Tl{\ell+1}:= \texttt{D\"orfler\_marking}(\theta, \sigma_\ell(K): K \in \Tl \ell)$}
		}
\end{algof}

%% file: algSAfem.tex
	\begin{algof}
	\DontPrintSemicolon
	\For{$\ell=0,1,\ldots$}{
		\acompute{$\eta_\ell(K)$, $\mu(K)$ for all $K \in \Tl \ell$\vspace{0.5ex}}
		\eIf(\tcp*[f]{Case (A) \hspace{1.5cm}}){$\estV{\ell}^2:=\mu^2(\Tl \ell) \leq \kappa\estJ{\ell}^2\equiv \kappa \eta_\ell^2(\Tl \ell)\eta^2(\Tl \ell, \Tl \ell)$}
		{
			{$\Tl{\ell+1}:= \texttt{D\"orfler\_marking}(\theta_A, \eta_\ell(K):K\in\Tl \ell)$\vspace{0.5ex}}
		}
		(\tcp*[f]{Case (B) \hspace{1.5cm}}){
		{
		$\Tl {\ell+1}:=\Tl \ell \oplus \texttt{appx}(\rho_B\mu^2_\ell, \mu(K):K\in\Tl 0)$}
		}
		}
\end{algof}

%% file: prelim.tex
\section{Axioms and results}\label{sec:axioms}
The axioms concern general conditions of the estimators $\eta$ and $\mu$,
which play different roles in the adaptive algorithm,  and are based on the set $\TO$ 
of admissible triangulations.

\subsection{Partitions and admissible triangulations}
\begin{figure}
	\begin{center}
	\includegraphics{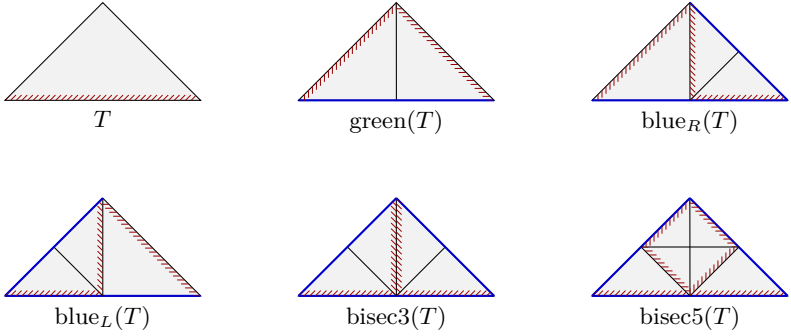}
	\end{center}
		\caption{Possible refinements of a triangle depending on the set of marked edges by NVB. Refinement edges are marked red, while marked edges are colored in blue.}
	\label{fig:refine}
\end{figure}

Let $\Tl 0$ be a regular triangulation of the domain $\Omega$ into (tagged) $n$-simplices in $\mathbb R^n$. 
Any refinement $\P$ from $\Tl 0$ by the newest vertex bisection (NVB) of \cref{fig:refine} is called partition, written $\P\in\PP{\Tl 0}=:\POO$. A partition $\P \in \POO$, which is a regular triangulation in the sense of Ciarlet, is called admissible, written $\P\in \TT{\Tl 0}=:\TO$.

The input of the underlying refinement procedure 
$\Tl{\text{out}}:=\textsc{Refine}(\Tl{\text{in}},\M)$ is an admissible triangulation  $\Tl{\text{in}}\in \TO$ and some subset
$\M\subseteq\Tl{\text{in}}$  thereof; the output $\Tl{\text{out}}$ is an admissible triangulation and a one-level refinement of $\Tl{\text{in}}$ with 
$\M\subset\Tl{\text{in}}\setminus\Tl{\text{out}}$ of quasi-minimal cardinality. Conversely, the procedure $\textsc{Refine}$ specifies 
the NVB with completion (to avoid hanging nodes etc.) and more details may be found in \cite{Stev08}. NVB is assumed throughout this paper. In particular, given $\T,\T'\in\TO$, their overlay $\T\oplus\T' \in \TT{\T} \cap \TT{\T'}$ is the smallest common 
refinement of $\T$ and $\T'$.

\subsection{Estimators and distance}\label{ssec:estDist}
The axioms are defined in terms of $\eta$ and $\mu$ plus a global distance $\delta$. For any admissible triangulation $\T\in\TO$ and any element domain $K\in\T$  let $\eta(\T,K)$ and $\mu(K)$ be a non-negative real number with squares $\eta^2(\T,K)$ and $\mu^2(K)$ and  their sums 
\begin{align}\label{eq:def_estim}
\eta^2(\T,\M)&:=\sum_{K\in\M}\eta^2(\T,K), \qquad \mu^2(\M):=\sum_{K\in \M}\mu^2( K) \qquad 
\text{for any }  \M\subseteq \T.
\end{align}
The distance $\delta(\T,\hT)$ of $\T\in\TO$ and its refinement $\hT\in\TO(\T)$ is a non-negative real.  
The estimators are utilized in the adaptive algorithm and are linked with the distance function in the axioms
below. The output of the adaptive algorithm is a sequence $\T_0,\T_1,\T_2,\dots$ of successive
refinements that start with $\T_0$ and give rise to the abbreviations (with a subindex $\ell$
to refer to the triangulation as part of the output of \safem)
\[
\eta_\ell(K):=\eta(\T_\ell,K)\quad\text{for }K\in\T_\ell\quad\mbox{and}\quad \eta_\ell:=\eta(\T_\ell,\T_\ell).
\]
The sum $\sigma^2:=\eta^2+\mu^2$ and their local variants are frequently utilized throughout this paper with $\sigma^2_\ell:=\eta^2_\ell+\mu^2_\ell$ for $\mu_\ell^2:=\mu^2(\Tl \ell):=\sum_{K\in \Tl \ell} \mu^2(K)$.

\subsection{Adaptive algorithm} 
In some more details, \safem calls \textsc{Select} and \textsc{Refine} to realize the D\"orfler marking in Case (A) from the introduction; more details on \texttt{appx} in Case (B) follow in Subsection \ref{sec:optApprox}.
\medskip

\centerline{\begin{minipage}{0.9\textwidth}
	\textbf{SAFEM$(\theta_A, \kappa, \rho_B, \Tl 0)$}
	\\
	\input{algAfem}

\end{minipage}}
The selection of $\Ml \ell$ with almost minimal cardinality means that $\abs{\Ml \ell} \lesssim \abs{\Ml \ell ^\star}$, where $\Ml \ell ^\star$ denotes some set of minimal cardinality with \cref{sAfem:eq:bulkA}. The point is that this can be realised in linear CPU time \cite{Stev07}.

\subsection{Axioms}
The universal positive constants $\Lref$, $\Lstab$, $\Lred$, $\LdRel$,  $\Lqo$, $\LsubAdd$, and 
 $\LhdRel\geq 0$ as well as  $0<\rho_2<1$ in the axioms (A1)-(A4), (B2), and (QM) 
below solely depend on $\TO$ (whence merely on $\Tl 0$); the parameters 
 $s>0$  and  $\LoptData$ in (B1)  also depend on the algorithm \texttt{appx} and the optimal data 
 approximation rate.
 
The axioms 
(A1)-(A3) and (B2) concern an arbitrary triangulation $ \T \in \TO$ and any refinement 
$\hT \in \TT{\T}$ of it, while (A4) solely concerns the outcome of \safem. Recall the sum conventions for $\eta(\T, \M)$ and $\mu(\T)$ in Subsection~\ref{ssec:estDist}.

\begin{axioms}
\item[(A1)]Stability. $\forall \T \in \TO \, \forall \hT \in \TT{\T}$
	\begin{equation}\label{eq:A1}\tag{A1}
		\abs{\eta(\hT,\T \cap \hT) - \eta(\T,\T \cap \hT)} \leq \Lstab \delta(\T,\hT).
	\end{equation}
\item[(A2)] Reduction. $\forall \T \in \TO\, \forall \hT \in \TT{\T}$
	\begin{equation}\label{eq:A2}\tag{A2}
		\eta(\hT, \hT \setminus \T) \leq \rho_{2}\eta(\T, \T \setminus \hT) + \Lred \delta(\T,\hT).
	\end{equation}
\item[(A3)] Discrete Reliability.
	$\forall \T \in \TO \,\forall \hT \in \TT{\T}\, \exists \mathcal R(\T, \hT) \subseteq \T$ with $\T \setminus \hT \subseteq \mathcal R(\T, \hT)$,
	\begin{align}\label{eq:A3}\tag{A3}
		\begin{aligned}
		\abs{ \mathcal R(\T, \hT)} &\leq \Lref \abs{\T \setminus \hT} \text{ and }\\
		\delta^2(\T, \hT) &\leq \LdRel \left( \eta^2(\T,\mathcal R(\T,\hT)) + \mu^2(\T) \right) + 
		\LhdRel \eta^2(\hT). 
	\end{aligned}
	\end{align}
\item[(A4)] Quasiorthogonality of discrete solutions.
	$\forall \ell \in \NO$
		\begin{equation}\label{eq:A4}
			\sum_{k=\ell}^{\infty} \delta^2(\Tl k, \Tl {k+1})
			\leq \Lqo \sigma_{\ell}^2.\tag{A4}
		\end{equation}
		\end{axioms}
\begin{axioms}
\item[(B1)] Rate $s$ data approximation.
	 $\forall \Tol>0$, $\Tl \Tol := \texttt{appx}(\Tol,\mu(K):K\in\Tl 0) \in \TO$  
	 satisfies
	\begin{equation}\label{eq:B}\tag{B1}
		\abs{\Tl \Tol}-\abs{\Tl 0} \leq \LoptData \Tol^{-1/(2s)} \quad \text{and} \quad  \mu^2(\Tl \Tol) \leq \Tol.
		\end{equation}
	\item[(B2)] Quasimonotonicity of $\mu$. $\forall \T \in \TO$ $\forall \hat \T \in \TT{\T}$\quad
 $\mu(\hat \T) \leq \LsubAdd \mu(\T)$.
\end{axioms}

Theorem~\ref{thm:qmono}  below asserts that the aforementioned axioms imply
quasimonotonicity of $\sigma$ for small values of  $\LhdRel $, while this axiom 
(QM) stands on its own in the example of Section~\ref{s:appllsfem}.
 
\begin{axioms}
\item[(QM)] Quasimonotonicity of $\sigma$. $\forall \T \in \TO$ $\forall \hT \in \TT{\T}$\quad
 $ \sigma(\hT) \leq \Lqm \sigma(\T)$.
\end{axioms}

\subsection{Optimal convergence rates}
The axioms\hspace{0.5mm}(A1)-(A4),%
\hspace{0.5mm}(B1)-(B2),\hspace{0.5mm}and\hspace{0.5mm}(QM) ensure 
quasioptimality of \safem for sufficiently small
parameters $\theta_A$ and $\kappa$ as stated in \cref{thm:safem} below. Recall that $\sigma^2:=\eta^2+\mu^2$ and set
\[
\sigma^2(\T)\equiv \sigma(\T)^2:=\sigma^2(\T,\T):=\sum_{K\in\T} \sigma^2(\T,K)
\text{ for $\T\in\TO$ and }
\sigma_\ell:=\sigma(\T_\ell).
\]
For any $N\in\NO$, the comparison with the optimal rates concern the optimal 
value 
\[
\min\sigma(\TT{N}):=\min\{\sigma(\T):\T\in\TT{N}\}
\] 
of all admissible triangulations 
\[
\TO(N):=\{\T\in\TO: |\T|\le |\T_0|+N\}
\]
of cardinality $|\T| \le |\T_0|+N$  
with at most $N$ extra cells.

\begin{theorem}[Quasioptimality]\label{thm:safem}
Suppose (A1)-(A4) and (B1)-(B2).	
\begin{inparaenum}[(a)]
\item
The strict inequality $(\Lstab^2+\Lred^2)\LhdRel < 1$ implies (QM) 
with $\Lqm$ depending on $\Lstab$, $\Lred$, $\LdRel$, $\LhdRel$, and $\LsubAdd$. 
\item The axiom (QM) leads to the existence of some $\kappa_0>0$,   which is $+\infty$ if
$\LsubAdd=1$, such that any choice of $\kappa$, $\theta_A$, and $\rho_B$ with 
\[
0<\kappa<\kappa_1:=\min\left\lbrace\kappa_0, \Lstab^{-2}\LdRel^{-1} \right\rbrace,\quad 
0<\theta_A<\theta_0:= (1-\kappa \Lstab^2\LdRel)  /(1+\Lstab^2\LdRel),
\] and  $0<\rho_B<1$ implies the following. 
The output  $(\Tl \ell)_{\ell \in \mathbb N_0}$ and $(\sigma_\ell)_{\ell \in \mathbb N_0}$ of \safem  satisfy the equivalence 
\begin{align}\label{eq:optim}
		\LoptData^s+ \sup_{\ell \in \NO}\left(1+\abs{\Tl \ell} - \abs{\Tl 0}\right)^{s}  \sigma_\ell 
		\approx \LoptData^s+\sup_{N\in \NO} (1+N)^s \min \sigma(\TT{N}).
\end{align}
\end{inparaenum}
\end{theorem}

In particular, the left-hand side of the equivalence \cref{eq:optim} is smaller than infinity if the  
right-hand is and vice versa. The quotient 
is bounded below and from above by the equivalence constants, which depend on 
$\Lref$, $\Lstab$, $\Lred$, $\LdRel$, $\LhdRel$, $\Lqo$, $\LsubAdd$, $\rho_B$, $\rho_2$,  
$\theta_A$, $\kappa$, and $s$ but not on $\Lambda_5$.

The (possibly unknown) parameter $s$ is not utilized in \safem.
The axioms (B1)-(B2) specify {\em sufficient} conditions for optimal 
convergence, where the parameter 
$s>0$ is arbitrary and may refer to a related  nonlinear approximation class.

%% file: algAfem.tex
	\begin{algof}
	\DontPrintSemicolon
	\KwIn{Initial coarse triangulation $\Tl 0$,
	$0<\theta_A< 1$, $0<\rho_B<1$,
	$0<\kappa$}
	\For{$\ell=0,1,\ldots$}{
		\acompute{refinement indicators $\estJ
		\ell^2(K)$ and  $\mu^2(K)$ for all $K \in \Tl \ell$ 
	\vspace{0.5ex}}
	\eIf(\tcp*[f]{Case (A)\hspace{5cm}}){$\estV{\ell}^2 \leq \kappa\estJ{\ell}^2$ }
		{
		\aselect{a subset $\Ml \ell
		\subseteq \Tl \ell$ of element domains 	of (almost) minimal cardinality with 
		\vspace{-1.5ex}
		\begin{align}\label{sAfem:eq:bulkA}
			\theta_A \estJ \ell^2 \leq  \estJ \ell
			^2\left(\Ml \ell\right):= \sum_{K \in \Ml \ell} \eta_\ell^2(K)
		\end{align}}
		\vspace{-5ex}
		\acompute{$\Tl {\ell+1}:= \textsc{Refine}(\Tl \ell,\Ml \ell)$ \vspace{0.5ex}}
		}(\tcp*[f]{Case (B)\hspace{5cm}}){
			\arun{$\T=\texttt{appx}(\Tol,\mu(K):K\in\Tl 0)$
			 with $ \Tol= \rho_B \estV{\ell}^2$}
		\acompute{$\Tl{\ell+1}:= \Tl \ell \oplus \T$}
		}}
		\KwOut{$\Tl k$, $\eta_k$, $\mu_k$, $\sigma_k:=\sqrt{\eta_k^2 + \mu_k^2}$ for $k=0,1,\ldots$}
\end{algof}

%% file: remarks.tex
	\section{Remarks}\label{sec:remarks}
	\subsection{Weak form of (A4)}
	The axiom (A4) can be a weakened with some parameter $\varepsilon>0$, which vanishes in 
	(A4)$\equiv$(A4$_0$).
	\begin{axioms}
		\item[(A4$_\varepsilon$)] Quasiorthogonality with $\varepsilon>0$.
		$\exists 
		\varepsilon>0$ $\exists 0<\LqoE<\infty \, \forall \ell,m \in \NO$
			\begin{equation}\label{eq:A4e}
				\sum_{k=\ell}^{\ell+m} \delta^2_{k, k+1} \leq \LqoE \sigma_{\ell}^2  + \varepsilon \sum_{k=\ell}^{\ell+m} \sigma_{k}^2  \tag{A4$_\varepsilon$}.
			\end{equation}
		\end{axioms}

	The axiom (A4$_\varepsilon$) implies (A4$_{\varepsilon'}$) for all $0\leq\varepsilon<\varepsilon'$ 
	with the same constant $\LqoE=\LqoEs$, and (A4) is (A4$_0$), i.e.\  (A4$_\varepsilon$) for $\varepsilon=0$). Conversely, as $\varepsilon \searrow 0$ it may be expected 
	that $\LqoE \rightarrow \infty$. In the presence of (A1)-(A2), this is not the case.  In fact, (A1)-(A2) 
	and (A4$_\varepsilon$) imply (A4) for sufficiently small $\varepsilon>0$.

	\begin{theorem}[(A4$_\varepsilon$)$\Rightarrow$(A4)]\label{thm:A4eA4} Let $\theta_A$ be the parameter of \safem and $0<\rho_{12}<1$  the reduction factor for the total error estimator with constant $0<\LII<\infty$ in \cref{thm:contraction} below and let
		$0\leq\varepsilon< (1-\rho_{12})/\LII$.
		Then (A1)-(A2) and (A4$_\varepsilon$) imply (A4)
		with $\Lqo:=\LqoE+\varepsilon (1+\LII \LqoE)/(1-\rho_{12}-\varepsilon \LII)$.
	\end{theorem}

	This has first been observed  in \cite{CFP14} for \cafem and is proved in 
	Subsection~\ref{ssec:conv} for completeness and applied below in \cref{thm:mfemA3}. 

	\subsection{Quasimonotonicity}\label{ssec:remQM} 
	The axiom (B2) explicitly ensures the 
	quasimonotonicity of $\mu$ and (QM) follows 
	with  $\Lqm:=\sqrt{\LsubAdd^2+\Lqmeta^2}$
	from the subsequent theorem: $\LhdRel < 1/(\Lstab^2+\Lred^2)$ is sufficient for  (QM). 

	\begin{theorem}[Quasimonotonicity] \label{thm:qmono}
			Suppose (A1)-(A3) and  $\widehat M:=(\Lstab^2+\Lred^2)\LhdRel$ $<1$. Set $M:=(\Lstab^2+\Lred^2)\LdRel$ and
			\begin{align*}
			\Lqmeta:&=\frac{1+M(1-\whM)+\whM +2 \sqrt{M (1-\whM)+\whM}}{(1-\whM)^2}.
	\end{align*}
	Then, any $\T \in \TO$ and $\hT \in \TT{\T}$ satisfy 
	\begin{align}\label{eq:qmono}
			\eta(\hT) &\leq \Lqmeta \sigma(\T). 
	\end{align}
	\end{theorem}
\textit{Proof.}
		Given $	\lambda:=(\sqrt{M+\whM-M \whM}-\whM)/(M+\whM)<1/\whM-1$,
	recall the following implication of the axioms \eqref{eq:A1}-\eqref{eq:A3}, namely
	\begin{align*}
		\eta^2(\hT,\hT \cap \T)&\leq (1+1/\lambda) \eta^2(\T,\hT \cap \T) + (1+\lambda)\Lstab^2 \delta^2(\T,\hT),
		\\
		\eta^2(\hT,\hT \setminus\T)) & \leq (1+1/\lambda)\rho_{2}^2 \eta^2(\T, \T \setminus \hT) +(1+\lambda) \Lred^2\delta^2(\T,\hT), 
		\\
	\delta^2(\T,\hT) &\leq \LdRel \sigma^2(\T) + \LhdRel \eta^2(\hT).
	\end{align*}
	Those inequalities plus the split 
	$\eta^2(\hT)= \eta^2(\hT,\hT \cap \T) + \eta^2(\hT,\hT \setminus \T))$
	 verify
	\begin{align*}
		\eta^2(\hT) &\leq  (1+1/\lambda)\eta^2(\T) + (1+\lambda)(\Lstab^2+\Lred^2)\left(\LdRel \sigma^2(\T)+\LhdRel \eta^2(\hT)\right).\qquad \endproof
		\end{align*}

\subsection{Optimal data approximation with \Approx}\label{sec:optApprox}
Case (B) of \safem runs a data approximation algorithm \texttt{appx}($\Tol, \mu(K):K \in \Tl 0$) 
with output in $\TO$. The  data approximation algorithm   \Approx  \cite{BDD04,BD04} is based on the refinement of partitions and has been  established for separate marking algorithms   in \cite{CR09,safem2015} and is one possible realisation of \texttt{appx} in \safem. 

Let $\hat \P$ be some NVB refinement of $\P \in \POO$.  Let $K\in \P$ and $\hat \P \in \PP{\P}$, then the refinement of $K$ in $\hat \P$ 
is the set $\hat \P(K):= \lbrace T \in \hat \P \, \vert \, T\subseteq K \rbrace$  in 
the following.

\begin{axioms}
\item[(SA)] Sub-additivity. 
$\exists \LsubAdd<\infty \, \forall \P \in \POO \, \forall \hat \P \in \PP{\P} \, \forall \M \subseteq \P$ 
	\begin{align}\label{eq:B1}\tag{SA}
		\mu^2(\hat \P(\M)):=\sum_{K \in \M} \sum_{T \in \hat \P(K)} \mu^2(T) 
				\leq \LsubAdd \mu^2(\M).	
 	\end{align}
\end{axioms}
Note, that the notation of the data approximation term $\mu$ is a straight forward extension of its definition in \eqref{eq:def_estim} for admissible triangulations to partitions.

The algorithm \Approx is outlined in the following with input tolerance 
$\Tol':= \Tol/\LsubAdd=\rho_B \mu_\ell / \LsubAdd$ and the values $\mu(K)$ on the coarse triangulation $\Tl 0$.
\medskip

\textbf{\Approx($\Tol', \mu(K):K\in\Tl 0$)} \\
\input{algTSA}

\begin{remark}
	\begin{inparaenum}[(a)]
	\item Algorithm \Approx is based on a modified error functional $\tilde \mu$ initiated by
$
	\tilde \mu (K):= \mu(K) \text{ for all } K \in \Tl 0. 
$
Given $\tilde \mu(K)$ for a triangle $K=K_1\cup K_2$ bisected  
into sub-triangles $K_1$ and $K_2$, let 
\begin{align}	
	\tilde \mu (K_j):= \frac{\tilde \mu (K)( \mu(K_1)+ \mu(K_2))}{\mu(K) + \tilde \mu (K)}
	\quad\text{ for } j=1,2 .\label{eq:tmu2}
\end{align}
\item 
Notice that the partitions $\P$ in the while-loop in  \Approx 
are not regular in general and the final completion step may be 
realized with  successive calls of 
\textsc{Refine}.
\\
\item The implementation of \Approx may store the partition $\P$ and the values
$\tilde \mu(K)$ for all element domains $K\in\P$ at the end of the while loop to keep the successive calls of \Approx for various decreasing tolerances $\Tol'$ efficient.
\end{inparaenum}
\end{remark}

\medskip

\begin{theorem}[\cite{BD04,BDD04}]\label{thm:optimalityofApprox}
(SA) in  \Approx implies  (B1)-(B2) with rate-s-optimality in the sense that
\begin{equation}\label{eqdefM(s,mu)}
M(s,\mu):=\sup_{N\in \NO} (1+N)^s \min \mu(\TT{N}) \approx \LoptData^s
\end{equation}
holds for all $s>0$ (and $M(s,\mu)<\infty$ if and only if $ \LoptData<\infty$).
\end{theorem}

\textit{Proof.}
This follows from near optimality proven in  \cite[Theorem 6.1]{BD04} and 
\cite[Lemma 4.4]{BDD04}. 
\qquad\endproof

\subsection{Collective D\"orfler marking is optimal for $\volO{\ell}$}\label{ssec:Doerfler}
Given $f\in \LO$ in the polyhedral domain $\Omega \subseteq \mathbb R^n$ partitioned
into the regular triangulation $\Tl 0$, set $\eta(\Tl \ell,K):=\abs{K}^{2/n} \vert f \vert_{L^2(K)}$ for all $K \in \Tl \ell$.
Let $\eta_\ell=\eta(\Tl \ell, \Tl \ell)$. Then, (A1)-(A4) are satisfied with appropriate weight functions $h_{\T}$ (resp.\ $h_{\hT}$) of mesh-sizes in $\T$ (resp.\ $\hT$)
\begin{align*}
	\delta(\T,\hT):= \norm{(h_{\T}-h_{\hT})f}{L^2(\Omega)}.
\end{align*}
Hence \cafem with collective D\"orfler marking implies optimal data approximation
for this particular data error term with a mesh-size weight $h_\T$.
This is in agreement with the well-established fact that  first-order conforming and nonconforming finite element methods do not need  a data reduction with \safem. 

%% file: algTSA.tex
\begin{algof}
	\DontPrintSemicolon
	\acompute{$\tilde \mu^2(T)=\mu^2(T)$  for all $T \in \Tl 0$ and set $\mu^2(\Tl 0):=\sum_{T \in \Tl 0}\mu^2(T)$}
	\aset{${\P} = \Tl 0$}
		\While{$\mu^2({\P})> \Tol'$}{
		\acompute{$\tilde \mu^2(T)$ for all $T \in {P}$, set ${\tilde \mu^2}_{\max}:=\max_{T \in {\P}} \tilde \mu^2(T)$}
		\aselect{a subset $\M := \left\{ T \in {\P} \mid \tilde \mu^2(T) = {\tilde \mu^2}_{\max} \right\}\subseteq {\P}$}
		\acompute{$\overline{\P}:=\texttt{bisec}(\P,\M)$}
		}
		\acompute{$\Tl \Tol:= \texttt{completion}({\P})\in \TO$}
\end{algof}

%% file: optimality.tex
\section{Proofs} \label{sec:optimality}
The abbreviation $\delta_{\ell,\ell+1}:=\delta(\Tl \ell, \Tl {\ell+1})$ applies throughout this section.

\subsection{Estimator reduction}
The constant $\LsubAdd \geq 1$ in the following theorem leads to $\kappa_0$ set to $+\infty$ for $\LsubAdd=1$; 
$\kappa_0=\infty$ and $\LsubAdd=1$ hold in all the examples of this paper.

\begin{theorem}[(A12) reduction] \label{thm:contraction}
Suppose \eqref{eq:A1}-\eqref{eq:A2} and parameters $0<\theta_A\leq 1$, 
$0<\kappa $,  and $0<\rho_B<1/ \LsubAdd$ 
from \safem. Any choice of $\gamma$ and $\lambda$ with 
	\begin{align}
	&0<\gamma<\rho_2^{-2}-1 \text{ and }0<\lambda 
	<\min \left\{\left( 1-(1+\gamma)\rho_2^2\right)\frac{\theta_A}{1-\theta_A}, 
	\kappa (1-\rho_B)\right\} \label{eq:gamma}\\
		\intertext{lead to constants }
&0<\LII:=(1+1/\lambda)\Lstab^2+(1+1/\gamma)\Lred^2<\infty \label{eq:LII},\\
&0<\rho_A:=(1+\lambda)(1-\theta_A) +(1+\gamma)\rho_{2}^2\theta_A<1, \label{eq:rhoA}\\
&0<\kappa_0:=(1-\rho_A)/(\LsubAdd-1) \text{ (with $\kappa_0:=+\infty$ if $\LsubAdd=1$)},\\
&0<\rho_{12}:=\max \left\lbrace \rho_{A} + \kappa\LsubAdd, 1+\lambda +\kappa
 \rho_B\right\rbrace/(1+\kappa)\le 1 .\label{eq:rho}
	\end{align}
Moreover, $0<\kappa  <\kappa_0$ implies $\rho_{12}<1$ and
	\begin{equation}
		\sigma_{\ell+1}^2 \leq \rho_{12} \sigma_{\ell}^2  + \LII \delta^2_{\ell, \ell+1} 
		\quad \text{for all } \ell\in\NO
		\tag{A12} \label{eq:A12}
	\end{equation}
for 	the output $\sigma_\ell^2$ of \safem.
\end{theorem}

\textit{Proof}
		For $\gamma$ and $\lambda$ as in \eqref{eq:gamma}, the axioms \eqref{eq:A1}-\eqref{eq:A2} imply 
		\begin{align*}
			\eta_{\ell+1}^2(\Tl {\ell+1} \cap \Tl \ell)  &\leq (1+\lambda) \eta_\ell^2(\Tl {\ell+1} \cap \Tl \ell) + (1+1/\lambda) \Lstab^2 \delta^2_{\ell,\ell+1}, \\ 
			\eta_{\ell+1}^2(\Tl {\ell+1} \setminus \Tl \ell) &\leq (1+\gamma)\rho_{2}^2\eta_{\ell}^2(\Tl \ell \setminus \Tl {\ell+1}) + (1+1/\gamma)\Lred^2 \delta^2_{\ell,\ell+1}. 
		\end{align*}
		The sum of those two inequalities leads to
		\begin{align}\label{eq:A1+A2}
			\eta_{\ell+1}^2 \leq (1+\lambda) \eta_\ell^2  + ((1+\gamma)\rho_{2}^2-(1+\lambda))\eta_{\ell}^2(\Tl \ell \setminus \Tl {\ell+1}) +  \LII \delta^2_{\ell,\ell+1}.
		\end{align}
		The restrictions on $\lambda$ and $\gamma$ ensure $(1+\gamma)\rho_{2}^2<1<1+\lambda $. 
		Thus, in general,
		\begin{align*}
			\eta_{\ell+1}^2\leq (1+\lambda)\eta_\ell^2 + \LII \delta^2_{\ell,\ell+1}.
		\end{align*}
		In Case (A) on the level $\ell$, when D\"orfler's marking ensures $ \theta_A \eta_\ell^2 \leq \eta_\ell^2(\Tl \ell \setminus \Tl {\ell+1})$, this and \eqref{eq:A1+A2} leads to an improvement of the last estimate, namely
		\begin{align*}
			\eta_{\ell+1}^2 &\leq  \left( (1+\lambda)(1-\theta_A)+ (1+\gamma)\rho_2^2\theta_A \right) \eta_{\ell}^2 + \LII \delta^2_{\ell,\ell+1}=\rho_A  \eta_{\ell}^2 + \LII \delta^2_{\ell,\ell+1}.
		\end{align*}
		The restrictions on $\lambda$ and $\gamma$ reveal $\rho_A <1$. Altogether, let
		\begin{align} \label{eq:LestRed}
			\LestRed&:= \begin{cases}
			\rho_A & \text{ in Case (A) on level $\ell$},\\
			1+\lambda & \text{ in Case (B) on level $\ell$}.
			\end{cases}
		\end{align}
		Then, the output of \safem  satisfies
		\begin{align}
			\estJ{\ell+1}^2 &\leq \LestRed \estJ{\ell}^2 + \LII \delta^2_{\ell,\ell+1}
			\quad \text{for all } \ell\in \NO. \label{eq:estRed}
		\end{align}
%%%%
		In Case (A) on any level $\ell$ with $\LestRed=\rho_{A}$ from \eqref{eq:rhoA} and $\LII$ from \eqref{eq:LII}, it also holds $\estV{\ell+1}^2 \leq \LsubAdd\estV{\ell}^2$, and $\estV{\ell}^2 \leq \kappa \estJ{\ell}^2$.
	Since $\alpha:=(\LsubAdd-\rho_A)/(\kappa+1) >0$, this and \cref{eq:estRed} lead to
	\begin{align*}
		\sigma_{\ell+1}^2 &\leq (\rho_{A} +\alpha \kappa )\eta_\ell^2 
		+ (\LsubAdd-\alpha) \mu_\ell^2+ \LII \delta^2_{\ell, \ell+1}
		=\frac{\rho_A+\kappa\LsubAdd}{1+\kappa}\sigma_\ell^2+ \LII \delta^2_{\ell, \ell+1}.
	\end{align*}
	In Case (B) on the level $\ell$ with $\LestRed=1+\lambda$, it  holds $\estV{\ell+1}^2 \leq 
	 \rho_B \estV{\ell}^2$, and $\kappa \estJ{\ell}^2<\estV{\ell}^2$.
		Since $\beta:=\kappa(1+\lambda-\rho_B )/(1+\kappa)>0$, 
		this and \cref{eq:estRed} lead to
		\[
			\sigma_{\ell+1}^2 < (1+\lambda-\beta)\estJ{\ell}^2 
			+ (\rho_B + \beta/\kappa) \estV{\ell}^2 + \LII \delta^2_{\ell,\ell+1}
			= \frac{1+\kappa\rho_B +\lambda}{1+\kappa}\sigma_{\ell}^2  
			+ \LII \delta^2_{\ell,\ell+1}.
	\]
	This proves the total error estimator reduction \cref{eq:A12} with $\rho_{12}$ from \cref{eq:rho}. \qquad \endproof

\subsection{Convergence}\label{ssec:conv}
The plain convergence follows from the estimator reduction (A12) plus quasiorthogonality (A4).
\begin{theorem}\label{th:conv}
Suppose $0<\theta_A\leq 1$, $0<\kappa$, $0<\rho_B<1$, suppose (A4) and 
(A12) with constants  $0<\rho_{12}<1$ and $0<\LII<\infty$.
Then $\La:=(1+\LII \Lqo)/(1-\rho_{12})$,  $q:=\Lambda/(1+\Lambda)<1$,
and the output of \safem satisfy the following assertions (a)-(c).
	\begin{alphenum}
		\item (Plain convergence) $ \forall \ell, m \in \NO$ \quad
				$\displaystyle \sum_{k=\ell}^{\ell+m} \sigma_k^2 \leq \La \sigma_\ell^2.$
		\item (R-linear convergence on each level)  $\forall \ell,m \in \NO$ 
		\quad	$\displaystyle	\sigma_{\ell+m}^2 \leq \frac{q^m}{1-q} \sigma_\ell^2.$
	\item (Reciprocal sum) \label{lem:geomRow} $\forall s>0$ $\forall \ell\in \mathbb N$ 
		\quad 	$\displaystyle	\sum_{k=0}^{\ell-1} \sigma_k^{-1/s} \leq \frac{q^{1/(2s)} \sigma_\ell^{-1/s}}{(1-q)^{1/(2s)}(1-q^{1/(2s)})}.$
	\end{alphenum}
\end{theorem}
 \noindent\textit{Proof of (a).}
	For all $\ell$, $m \in \NO$, (A12) implies
	\begin{align}\label{eq:A12imply}
		\sum_{k=\ell}^{\ell+m} \sigma_k^2 & = \sigma_\ell^2 + \sum_{k=\ell+1}^{\ell+m} \sigma_k^2
		\leq \sigma_\ell^2 +  \rho_{12}\sum_{k=\ell}^{\ell+m} \sigma_{k}^2  + \LII \sum_{k=\ell}^{\ell+m}\delta^2_{k,k+1}.
		\end{align}
		This plus (A4) verify
		\begin{align*}
		 \left( 1- \rho_{12} \right)\sum_{k=\ell}^{\ell+m} \sigma_k^2& \leq \sigma_\ell^2 + \LII \Lqo \sigma_\ell^2 .
	\end{align*}
	This proves (a) with the asserted constant $\La$.
	\qquad \endproof

	\noindent\textit{Proof of \cref{thm:A4eA4}.} The same argument as in the proof of (a) before show that (A12) and (A4$_\varepsilon$) imply (A4) for small $\varepsilon$. In fact, \cref{eq:A12imply} and (A4$_\varepsilon$) show
	\begin{align*}
		(1-\rho_{12})\sum_{k=\ell}^{\ell+m} \sigma_k &\leq \sigma_\ell^2+\LII \left(\LqoE \sigma_\ell^2 + \varepsilon \sum_{k=\ell}^{\ell+m}\sigma_k \right).
		\intertext{In other words}
		\left( 1- \rho_{12} - \varepsilon\LII \right)\sum_{k=\ell}^{\ell+m} \sigma_k^2& \leq \left( 1 + \LII  \LqoE \right)\sigma_\ell^2.
	\end{align*}
	This plus \eqref{eq:A4e} lead to (A4) with $\Lqo:= \LqoE + \varepsilon (1+\LII\LqoE)/(1-\rho_{12}-\varepsilon\LII)$
	\begin{align*}
		\sum_{k=\ell}^{\ell+m}\delta^2_{k,k+1} &\leq \LqoE \sigma_\ell^2 + \varepsilon\sum_{k=\ell}^{\ell+m} \sigma_k^2 
		\leq \Lqo \sigma_\ell^2 .
		\qquad \endproof
	\end{align*}

	 \noindent\textit{Proof of Theorem \ref{th:conv}.b.}
	The assertion (a) implies the convergence of the series
	\begin{align*}
		\xi_{\ell+1}^2:=\sum_{k=\ell+1}^{\infty} \sigma_k^2 &\leq \La \sigma_\ell^2<\infty.
	\end{align*}
	The addition of $\La \xi_{\ell+1}^2$ to the previous inequality results in
	\begin{align}
		(\La +1)\xi_{\ell+1}^2 &\leq \La \xi_{\ell}^2, \text{ hence } \xi_{\ell+1}^2\leq q \xi_{\ell}^2. \label{eq:xiContr}
	\end{align}
	The successive application of the previous contraction \eqref{eq:xiContr} shows
		\begin{align*}
			\sigma_{\ell+m}^2\leq \xi_{\ell+m}^2\leq q^m \xi_{\ell}^2 = q^m \left( \sigma_\ell^2+\xi_{\ell+1}^2 \right)
			\leq q^m (1+\La) \sigma_\ell^2. \qquad \endproof
	\end{align*}

 \noindent\textit{Proof of Theorem \ref{th:conv}.c.} The R-linear convergence of (b) leads to
  	\begin{align*}
		\sigma_k^{-1/s} \leq \frac{q^{(\ell-k)/(2s)}}{(1-q)^{1/(2s)}}  \sigma_\ell^{-1/s}\qquad \text{for all } 0\leq k <\ell.
	\end{align*}
 	This proves 
	\begin{align*}
		\sum_{k=0}^{\ell-1} \sigma_k^{-1/s} &\leq \frac{\sigma_\ell^{-1/s}}{(1-q)^{1/(2s)}} \sum_{k=0}^{\ell-1}\left(q^{1/(2s)}\right)^{\ell-k} \leq \frac{\sigma_\ell^{-1/s} q^{1/(2s)}}{(1-q)^{1/(2s)}(1-q^{1/(2s)})}.\qquad \endproof
	\end{align*}
	
\begin{lemma}[Comparison]\label{lem:competition} 
Suppose \eqref{eq:A1}-\eqref{eq:A4}, (B1)-(B2) with $0<s<\infty$, (QM),
		$0<q<1$ from Theorem \ref{th:conv}.b, and let $0<\xi<1$ and $0<\nu<\infty$; let 
		\begin{equation} \label{eqdefM(s,sigma)}
			M:=M(s,\sigma):= \sup_{N\in \NO}(N+1)^s \min \sigma(\TT{N})<\infty,
		\end{equation}
		similar to the definition of $M(s,\mu)$ in \eqref{eqdefM(s,mu)}.
 		Then for any level $\ell \in \NO$ of \safem with a triangulation $\Tl \ell$,  
	there exists a refinement $\hTl \ell \in \TT{\Tl \ell}$ with (a)-(c).
	\begin{alphenum}
		\item $\displaystyle \sigma (\hTl \ell) \leq \xi \sigma_\ell$;
		\item $\displaystyle\sqrt{1-q}\xi \, \sigma_\ell \,  
		\abs{\Tl \ell \setminus \hTl \ell}^s \leq \Lqm M$;
		\item 
			$ \left( 1 -\xi^2(1+\nu+ (1+1/\nu)\Lstab^2\LhdRel)\right)\estJ{\ell}^2$  \\
			\phantom{xx} $\leq  
			\left(1+( 1+1/\nu)\Lstab^2\LdRel\right) \eta^2_\ell(\mathcal R(\Tl \ell, \hTl \ell))$
			\\
			\phantom{xxxx} $+ \left((1+\nu)\xi^2 + (1+1/\nu)\Lstab^2(\LdRel+\LhdRel\xi^2) \right)\mu^2_\ell .$

	\end{alphenum}
\end{lemma}

\textit{Proof.}
Two pathological situations are excluded in the beginning of the proof. 
First, if  $\sigma_\ell=0$, then $\hTl \ell=\Tl \ell$ satisfies the assumptions (a)-(c).
Second, Theorem \ref{th:conv} guarantees convergence of some sequence of triangulations and (QM) 
implies that this holds for uniform refinements as well. Hence there exists a refinement 
$\hTl \ell$ of $\Tl \ell$ with (a) and $\hTl \ell\cap\Tl \ell=\emptyset$. The latter implies (c)
even in case  $M\equiv M(s,\sigma)=\infty$ when (b) is obvious. 

Throughout the remaining parts of the proof, it is therefore assumed that $M<\infty$ and 
$\sigma_\ell>0$. Then (QM)  implies $0<\sigma_{0}\leq M<\infty$.

	\textit{1. Setup.}
	Let $N_\ell \in \NO$ be minimal with 
	\begin{align}\label{eq:Nell_min}
		(N_\ell+1)^{-s} \leq \frac{\xi \sqrt{1-q}}{\Lqm M}\, \sigma_\ell.
	\end{align}
	The quasimonotonicity (QM) followed by the definition of 
	 $M:=M(s,\sigma)<\infty$ in \eqref{eqdefM(s,sigma)} and $0<q<1,0<\xi<1$ lead to
	\begin{align*}
		\frac{\xi\sqrt{1-q}}{\Lqm} \sigma_{\ell} \leq \xi \sqrt{1-q} \, \sigma_0 \leq \xi \sqrt{1-q} M < M.	
	\end{align*}
	Hence, $(N_\ell+1)^{-s}<1$ and so $N_\ell\geq 1$.
		Since $N_\ell \in \mathbb N$ is minimal with \eqref{eq:Nell_min},
	\begin{align}\notag 
		0<(N_\ell+1)^{-s}\leq& \frac{\xi\sqrt{1-q}}{\Lqm M} \sigma_\ell < N_\ell^{-s}.\\
		\intertext{This implies}\label{eq:invNMinimal}
		N_\ell^s < &\frac{\Lqm M}{\xi\sqrt{1-q}} \sigma_\ell^{-1}.
	\end{align}
	
	\textit{2. Design of $\hTl \ell$.}
	The definition of $M<\infty$  yields the existence of some optimal $\tTl \ell \in \TT{N_\ell}$ with 
	\begin{align} \label{eq:lowM}
		\left( N_\ell+1 \right)^s \sigma(\tTl \ell) \leq M.
	\end{align}
	The overlay triangulation $\hTl \ell:= \Tl \ell \oplus \tTl \ell$ \cite{CKNS07,Stev07} satisfies
	\begin{align}\label{eq:overlay}
		\abs{\hTl \ell} + \abs{\Tl 0} \leq \abs{\Tl \ell} + \abs{\tTl \ell}.
	\end{align}

	\textit{3. Proof of (a).} The quasimonotonicity (QM) followed by \eqref{eq:lowM} and \eqref{eq:Nell_min} shows
	\begin{align*}
		\sigma(\hTl \ell) \leq {\Lqm}  \sigma(\tTl \ell) \leq {\Lqm M} 
		(N_\ell+1)^{-s} \leq \xi \sigma_\ell\sqrt{1-q}<\xi \sigma_\ell.\qquad \endproof
	\end{align*}

	\textit{4. Proof of (b).}
	The definition of $\tTl \ell$, the overlay estimate in \eqref{eq:overlay}, and the upper bound for $N_\ell$ in \eqref{eq:invNMinimal} lead to
	\begin{align*}
		\abs{\Tl  \ell\setminus \hTl \ell}\leq \abs{\hTl \ell}-\abs{\Tl \ell} \leq
		\abs{\tTl \ell}-\abs{\Tl 0}\leq N_\ell \leq \left( \frac{\Lqm M}{\xi \sigma_\ell\sqrt{1-q}} \right)^{1/s}.\quad \endproof
	\end{align*}
	
	\textit{5. Proof of (c).}
	For any $0<\nu< \infty, 0<\xi<1$,  (A1) and (A3) result in 
	\begin{align*}
		\eta^2_\ell(\Tl  \ell\cap \hTl \ell) &\leq (1+\nu)\eta^2(\hTl \ell,\hTl \ell \cap \Tl \ell) + (1+1/\nu)\Lstab^2 \delta^2(\Tl \ell,\hTl \ell)\\
		& \leq\left(1+\nu+ (1+1/\nu)\Lstab^2\LhdRel\right)\eta^2(\hTl \ell) \\
		& \quad + (1+1/\nu)\Lstab^2\LdRel \left( \eta^2_\ell(\mathcal R(\Tl \ell, \hTl \ell)) + \mu^2_\ell\right).%\\
	\end{align*}
	This, (a),  and $\Tl \ell \setminus \hTl \ell \subseteq \mathcal R(\Tl \ell, \hTl \ell)$ result in
	\begin{align*}
		\estJ{\ell}^2 &= \eta^2_\ell(\Tl \ell\cap \hTl \ell) + \eta^2_\ell(\Tl \ell \setminus \hTl \ell) \\
		&\leq \left(1+\nu+ (1+1/\nu)\Lstab^2\LhdRel\right)\xi^2 \sigma^2_\ell + \left(1+( 1+1/\nu)\Lstab^2\LdRel\right) \eta^2_\ell(\mathcal R(\Tl \ell, \hTl \ell)) \\
		&\quad +   ( 1+1/\nu)\Lstab^2\LdRel\mu^2_\ell. 
		\end{align*}
		Some rearrangements with $\sigma^2_\ell=\eta^2_\ell+\mu_\ell^2$ prove (c).
\qquad \endproof
\subsection{Proof of Theorem \ref{thm:safem}}\label{Subsectionname=ProofofTheorem}

\textit{Proof of Theorem \ref{thm:safem}.a.}
This is a consequence of Theorem \ref{thm:qmono} plus (B2).
\qquad \endproof

\textit{Proof of ``$\lesssim $'' in \eqref{eq:optim} of Theorem \ref{thm:safem}.b.}
	Since $\theta_A<\theta_0$ and the function 
	\begin{align*}
		f(\xi,\nu):=\frac{1 -\xi^2\left((1+\kappa)(1+\nu)+ (1+\kappa)(1+1/\nu)\Lstab^2\LhdRel\right) -\kappa (1+1/\nu)\Lstab^2\LdRel}{1+( 1+1/\nu)\Lstab^2\LdRel}
	\end{align*}
	is strictly smaller than $\theta_0=\lim_{\nu \to \infty} f(0,\nu)$, there exists 
	$\nu$, $\xi$ such that 
	\[\theta_A< f(\xi,\nu)< \theta_0. \]
	Given $\kappa_0$ from Theorem~\ref{thm:contraction} and assume 
	$\kappa<\kappa_1:=\min\left\lbrace\kappa_0,\Lstab^{-2}\LdRel^{-1}\right\rbrace$.

\textit{Case (A).} 
	Lemma \ref{lem:competition}.c and $\mu_\ell^2 \leq \kappa \eta_\ell^2$  prove 
	that $\mathcal R(\Tl \ell, \hTl \ell)$ satisfies 
		\begin{align*}
	 		&\left( 1 -(1+\kappa)\xi^2(1+\nu)- (1+\kappa)\xi^2(1+1/\nu)\Lstab^2\LhdRel -\kappa  (1+1/\nu)\Lstab^2\LdRel\right)\estJ{\ell}^2 \\
			&\qquad \leq  
		\left(1+( 1+1/\nu)\Lstab^2\LdRel\right) \eta^2_\ell( \mathcal R(\Tl \ell, \hTl \ell)).
		\end{align*}
		This reads $\theta_A\eta_\ell^2 \le  f(\xi,\nu)\eta_\ell^2\le \eta^2_\ell( \mathcal R(\Tl \ell, \hTl \ell))$ and 
 implies  that $\mathcal R(\Tl \ell, \hTl \ell)$ satisfies D\"orfler
		marking in Case (A).

		Let $\Ml \ell=:\Ml \ell^{(0)}$ be the set of marked elements in the D\"orfler 
		marking on level $\ell$, while $\Ml \ell^{\star}$ is the optimal set of marked elements. 
		Hence, there exists $0<\Lopt<\infty$ such that 
	\begin{align*}
		\abs{\Ml \ell} \leq \Lopt\abs{\Ml \ell^{\star}} \leq \Lopt\abs{\mathcal R(\Tl \ell, \hTl \ell)}.
	\end{align*}
	The control over  $\mathcal R(\Tl \ell, \hTl \ell)$ of (A3) in  
	Lemma \ref{lem:competition}.b results in
	\begin{align*}
		\abs{\mathcal R(\Tl \ell, \hTl \ell)} \leq  \Lref\abs{\Tl \ell \setminus \hTl \ell} \leq
		\Lref \left(\frac{\Lqm M }{\sqrt{1-q}\xi \sigma_\ell}\right)^{1/s}.
	\end{align*}
	Hence,  $\LoptA:=\Lopt\Lref \Lqm^{1/s}(\sqrt{1-q}\xi)^{-1/s}$ satisfies 
	\begin{align}\label{eq:optA}
		\abs{\Ml \ell^{(0)}}=\abs{\M _\ell} \leq \LoptA M^{1/s}\sigma_\ell^{-1/s}.
	\end{align}

\textit{Case (B).}
The output of \texttt{appx} with  input triangulation $\Tl 0$ and input tolerance 
	$\Tol:= \rho_B \mu_\ell^2$ on the level $\ell$ satisfies  (B1). Since 
	$\sigma_\ell^2 = \eta_\ell^2 + \mu_\ell^2 \leq (1+1/\kappa)\mu_\ell^2$ 
	in Case (B), this leads to 
\begin{align*}
		\abs{\Tl \Tol}-\abs{\Tl 0} \le  \LoptData (1+1/\kappa) \rho_B ^{-1/(2s)} \sigma_\ell^{-1/s} .
	\end{align*}
According  to \cite{CR09,safem2015} for $\Tl{\ell+1}= \Tl \ell \oplus \Tl \Tol$ 
there exists  a finite sequence 
$(\Ml \ell^{(k)})_{k=0,\ldots,K(\ell)}$ of sets of marked element domains that 
$\Tl{\ell}^{(0)}:=\Tl{\ell}$ and 
satisfies
\begin{align*}
		\Tl{\ell}^{(k+1)}= \text{\Refine}( \Tl \ell^{(k)}, \Ml \ell^{(k)})
		\quad\text{for all } k=0,\ldots,K(\ell)-1
\end{align*}
leads to $\T_{\ell+1}=\Tl{\ell}^{(K(\ell))}$. 
This observation and the estimate for the overlay with the sequence 
$(\Ml \ell^{(k)})_{k=0,\ldots,K(\ell)}$  \cite[Theorem 3.3]{CR09} show
\begin{align}\label{eq:optB}
		\sum_{k=0}^{K(\ell)} \vert\Ml \ell^{(k)}\vert \leq \abs{\Tl \Tol} - \abs{\Tl 0}
			\lesssim  \LoptData (1+1/\kappa) \rho_B ^{-1/(2s)}  \sigma_\ell^{-1/s}.
\end{align}
	The estimate from \cite[Theorem 3.3]{CR09} is for 2D only, however it is expected to hold in general.

\textit{Finish of the proof of ``$\lesssim$''.} 
	It is proven in \cite{CR09,safem2015} that the overhead control of \cite{BDD04,Stev08} holds in the sense that
	\begin{align}\label{eq:overheadClB}
		\abs{\Tl \ell}-\abs{\Tl 0 } \leq \Lcl
		 \sum_{j=0}^{\ell-1} \sum_{k=0}^{K(j)}\vert\Ml j^{(k)}\vert.
	\end{align}
	With \eqref{eq:optA}-\eqref{eq:optB} and Theorem \ref{th:conv}.c,  this proves
	\begin{align}\label{eq:optFin}
		\abs{\Tl \ell}-\abs{\Tl 0} &  
		\lesssim (\LoptData+M^{1/s}) \sigma_\ell^{-1/s}.
	\end{align}
	Finally, $1\leq \abs{\Tl \ell}- \abs{\Tl 0}$ implies $1+  \abs{\Tl \ell}- \abs{\Tl 0} \leq 2( \abs{\Tl \ell}- \abs{\Tl 0})$ while $\abs{\Tl \ell}=\abs{\Tl 0}$ implies $1\leq \sigma_\ell^{-1/s}(\LoptData+M^{1/s})$. Hence \eqref{eq:optFin} proves 
	$\sigma_\ell (1+ \abs{\Tl \ell}- \abs{\Tl 0})^s \lesssim \LoptData^s +M$ 
	and so ``$\lesssim$'' in the assertion of \cref{thm:safem}.
\qquad \endproof

\textit{Proof of ``$\gtrsim$''  in \eqref{eq:optim} of Theorem \ref{thm:safem}.b.}
Given $N\in\NO$ suppose that $\min\sigma(\TT N)$ is positive and so $\sigma_\ell>0$ for all 
$\ell\in\NO$ with $N_\ell:= |\T_{\ell}|-|\T_0|\le N$. This leads on the level $\ell$ in \safem to $N_{\ell+1}>N_\ell$ 
for it only stops with $\T_{\ell}=\T_{\ell+1}=\T_{\ell+2}=\dots$ when $\sigma_\ell=0$. Hence there exists some level $\ell$ with $N_\ell<N\le N_{\ell+1}$. This implies 
\begin{equation}\label{cceqN+1Sigmaell}
(N+1)^s\min\sigma(\TT N)\le (N_{\ell+1}+1)^s \sigma_\ell,
\end{equation}
which is evident in case $\min\sigma(\TT N)=0$. 

In Case (A) on the level $\ell$ of  \safem, there is a one-level refinement to create $\T_{\ell+1}$
(indicated in  \cref{fig:refine} for 2D), where each simplex in $\T_{\ell}$ creates a finite number $\le K(n)$ of children in a completion step. The constant $K(n)\ge 2$ depends only 
on the spatial dimension   $n$ \cite{GSStevcmam}. This leads to the bound 
$|\T_{\ell+1}|\le K(n)\, |\T_{\ell}| $ and then to
\[
(N_{\ell+1}+1)/(N_{\ell}+1)\le  K(n)+ (K(n)-1)(|\T_0|-1)\lesssim 1.
\]
In Case (B) on the level $\ell$ of  \safem, the refinement $\T_{\ell+1}:=\T_{\ell}\oplus \T_{\Tol}$
is controlled by $|\T_{\Tol}|-|\T_0|\le \LoptData \Tol^{-1/(2s)}\le  \LoptData\rho_B^{-1/(2s)}
\mu_{\ell}^{-1/s}$. Since  $\sigma_\ell^2\le (1+1/\kappa)\mu_{\ell}^2$ in Case (B), 
the overlay estimate
of \cite{CKNS07,Stev07} proves
\[
N_{\ell+1}-N_\ell\le |\T_{\Tol}|-|\T_0|\le \LoptData\rho_B^{-1/(2s)} (1+1/\kappa)^{1/(2s)}
\sigma_{\ell}^{-1/s}.
\]
This leads to the bound 
\[
2^{-s}(N_{\ell+1}+1)^s\le (N_\ell+1)^s + \rho_B^{-1/2} (1+1/\kappa)^{1/2} \LoptData.
\]
Consequently, in each of the  Cases (A) and (B), it follows 
\[
(N_{\ell+1}+1)^s\sigma_\ell\le   \left(K(n)+ (K(n)-1)(|\T_0|-1)\right)^s  (N_{\ell}+1)^{s}\sigma_\ell
+2^{s}\rho_B^{-s/2} (1+1/\kappa)^{s/2}\LoptData^s.
\]
With $S:= \sup_{\ell \in \NO}\left(N_{\ell}+1 \right)^{s}  \sigma_\ell$, this and   \eqref{cceqN+1Sigmaell} 
imply 
\begin{align*}
(N+1)^s\min\sigma(\TT N)\le  \left(K(n)+ (K(n)-1)(|\T_0|-1)\right)^s S + 
2^{s}\rho_B^{-s/2} (1+1/\kappa)^{s/2}\LoptData^s.
\end{align*}
Since this holds for any $N\in\NO$, the previous $N$-independent upper bound is greater than or equal to 
the supremum $M$ as well. This concludes the proof of ``$\gtrsim$''  in \eqref{eq:optim}.
\endproof

%% file: applications.tex
\section{Application to mixed FEM}
\label{s:applmfem}
The a~posteriori error analysis of mixed finite element schemes \cite{CC-MC-97,Alonso96} was completed in  \cite{ccdpas2015} with a reliable and efficient error control in  $H(\ddiv,\Omega)\times L^2(\Omega)$,
which is the natural functional analytical framework for the dual formulation of a Poisson model problem.

Given the right-hand side $f\in L^2(\Omega)$, the dual formulation of the Laplace equation 
on a 2D polygonal bounded simply-connected Lipschitz domain $\Omega$ 
seeks  $p\in H(\ddiv,\Omega)$ and $u\in L^2(\Omega)$ with 
\begin{align}\label{e:abstractmixed}
	\begin{aligned}
		a(p,q)+b(q,u)&=  0\quad\mbox{for all } q\in H(\ddiv,\Omega),  \\
		b(p,v)&=-F(v) := -\int_\Omega fv\, dx \quad\mbox{for all }v\in L^2(\Omega).
	\end{aligned}
\end{align}
Therein, the bilinear forms model the $L^2$ scalar product and the divergence term,
\begin{equation}
\label{e:bilinear}
a(p,q):=\int_\Omega p\cdot q\, dx\quad\text{and} \quad 
b(q,v):=\int_\Omega v\,\ddiv q\, dx.
\end{equation}
It is well established that  the weak solution 
$u\in V:= H^1_0(\Omega)$ to $-\Delta u=f$  in $\Omega$ 
specifies the flux $p:=\nabla u$; the two formulations are equivalent and allow for unique solutions.

Given an admissible triangulation $\T\in\TO$  let 
$(p_{RT},u_{RT})\in RT_0(\T)\times P_0(\T)$ solve the discrete problem 
\begin{align}\label{e:abstractmixedfem}
	\begin{aligned}
a(p_{RT},q_{RT})+b(q_{RT},u_{RT})&= 0\quad\mbox{for all } q_{RT}\in RT_0(\T),  \\
	b(p_{RT},v_{RT})&=-F(v_{RT})  \quad\mbox{for all }v_{RT}\in P_0(\T).
\end{aligned}
\end{align}
Given the unique discrete solution $(p_{RT},u_{RT})$ (resp. $(\hpRT,\huRT)$) with respect to the triangulation 
$\T\in\TO$ (resp.\ its refinement $\hT\in\TO(\T)$),  the estimators  of \cite{ccdpas2015} 
and the distance function read
\begin{align*}
\eta^2(\T,K)&:= |K|\, || p_{RT} ||_{L^2(K)}^2+ |K|^{1/2}\, \sum_{E\in \E(K)} ||[p_{RT}]_E\cdot\tau_E\|_{L^2(E)}^2,\\
\mu^2(K) & := ||f- f_K||_{L^2(K)}^2 \qquad\text{for any }K\in\T ,\\
\delta^2(\T,\hT) & := \| \hpRT-p_{RT}\|_{\Hdiv}^2. 
\end{align*} 
The standard 2D notation applies  to the triangle $K$ of area $|K|$ and its set $\E(K)$ of the three edges and the integral mean  $f_K:=\int f(x)\, dx/|K|$ of $f$. The jump $[\bullet]_E$ across  an interior edge $E=\partial T_+\cap\partial T_-$
with tangential normal vector $\tau_E$ and normal $\nu_E$ 
is the difference of the respective traces 
$ [q]_E:= q|_{T_+}-q|_{T_-}$ on $E$ from the two neighboring triangles $T_\pm$. Homogeneous 
Dirichlet boundary data translate into homogeneous jumps on the boundary:
 $ [q]_E:= q|_{T_+}$ for $E\subset\partial\Omega$ with neighboring triangle $T_+$. 
 
It is remarkable that, in the lowest-order case at hand, the Lagrange multiplier 
$u_{RT}$  does \textit{not} enter the estimators and hence the distance function acts on the flux approximations only. 

\begin{theorem}[(A1)-(A4)] \label{thm:mfemA3}
The estimators and distance functions satisfy (A1)-(A4) and (B2) for 
$\mathcal{R}(\T,\hT):=\T\setminus\hT$, $\Lref=1=\LsubAdd$, and $\LhdRel =0$. 
\end{theorem}

The estimator is reliable and efficient \cite{ccdpas2015} in that the exact (resp.\ discrete) solution $(p,u)$ (resp.\ $(p_{RT},u_{RT})$ with respect to $\T\in\TO$) satisfies
\[
\sigma(\T)\approx \| p-p_{RT}\|_{H(\ddiv,\Omega)}+\| u-u_{RT}\|_{L^2(\Omega)}.
\] 
Hence the optimal rates of the estimators is equivalent to the optimal rates of the errors in terms of nonlinear approximation 
classes with respect to the natural norms in $H(\ddiv)\times L^2$ of the mixed FEM.

\textit{Proof of \cref{thm:mfemA3}.}
It is straightforward to see that the estimators and distance function satisfy (A1)-(A2)
with $\rho_2:=2^{-1/4}$ and $\Lstab=\Lred\approx 1$ stemming from trace 
and inverse estimates.

The proof of (A3) requires an intermediate solution $\hpRT^*\in RT_0(\hT)$ with respect to the fine 
triangulation $\hT$ to the  above  Poisson model problem with a piecewise constant right-hand side 
$\Pi_0 f\in P_0(\T)$ with respect to the coarse triangulation $\T$.
Let $\mathcal{E}'\subseteq \E$ be the subset of all edges % with length $|E|$
such that at least one of the neighboring triangles
$K\in\T\setminus\hT$ with $E\in\E(K)$ is refined ($K\notin\hT$).
The  divergence-free Raviart-Thomas function $\hpRT^*-\pRT$ equals the rotated gradient
of some continuous and piecewise affine function and so gives rise to a stability result
\begin{equation}\label{stabilityintermediate}
\|   \hpRT^*-\pRT\|^2_{L^2(\Omega)}  \lesssim  \sum_{E\in \mathcal{E}'} |E|\, ||[\pRT]_E\cdot\tau_E ||_{L^2(E)}^2
\end{equation}
proved via a discrete Helmholtz decomposition (cf. e.g. \cite[Thm 5.6]{LCMHJX} for 
references and the arguments) for a simply connected domain $\Omega$.  

The discrete inf-sup condition (with respect to the finer mesh $\hT$) 
leads to some $\hqRT\in RT_0(\hT)$ and $\widehat{v_0}\in P_0(\hT)$ with norm
$ \| \hqRT\|_{H(\ddiv,\Omega)} +\| \widehat{v_0}\|_{L^2(\Omega)}\lesssim 1$ and 
\begin{align} \notag
	\text{LHS}:=&\,\| \hpRT-\pRT\|_{H(\ddiv,\Omega)} +\| \huRT-u_{RT} \|_{L^2(\Omega)}\\ 
	=& \,a(\hpRT-\pRT,\hqRT)+b(\hqRT, \huRT-u_{RT}  )+b( \hpRT-\pRT, \widehat{v_0}  ).
\intertext{The discrete equations \eqref{e:abstractmixedfem} on the fine level $\hT$ and $\ddiv p_{RT}=-\Pi_0 f$ show}
\text{LHS}=&-a(\pRT,\hqRT)-b( \hqRT, u_{RT} )- F( \widehat{v_0} -\Pi_0 \widehat{v_0}   ).
\label{eq:est:LHS}
\end{align}
Given $ \hqRT$ with bounded norm, let $ \qRT$ denote the mixed finite element solution to
a Poisson model problem with right-hand side $-\Pi_0 \ddiv  \hqRT \in P_0(\T)$. This leads to 
$\|\qRT\|_{H(\ddiv,\Omega)}\lesssim  1 $ and 
\begin{align*}
b( \hqRT, u_{RT} )=b( \qRT, u_{RT} )=-a(\pRT,\qRT).
\end{align*}
With $\|  \widehat{v_0}\|\lesssim 1$, the combination of the two previously displayed formulas shows
\[
\text{LHS}\lesssim
 \|\widehat{\Pi_0} f-\Pi_0 f\|_{L^2(\Omega)} + a(\pRT,\qRT-\hqRT).
\]
The Cauchy-Schwarz inequality leads to
\begin{align}\label{eq:biA:CSI}
a(\pRT,\qRT-\hqRT)&=a(\pRT-\hpRT^*,\qRT-\hqRT)+a(\hpRT^*,\qRT-\hqRT)\\
&\lesssim 
\|\pRT-\hpRT^*\|_{L^2(\Omega)}
 % \sqrt{ \sum_{E\in \mathcal{E}'} |E|\, ||[\pRT]_E\cdot\tau_E ||_{E}^2} 
 + a(\hpRT^*,\qRT-\hqRT).
\end{align}
Due to \eqref{stabilityintermediate} it remains to analyze the latter term.
The  test function equals \cite{Marini,ArnoldBrezzi}
\begin{equation}\label{aforementionedrepresentation}
 \qRT-\hqRT=\nabla_{NC} \widehat{v_{CR}} + \curl \widehat{\beta_C}+ 
 1/2\, \left((\Pi_0-1)\ddiv \hqRT \right)(\bullet -\midhT)
\end{equation}
for unique discrete functions $\widehat{v_{CR}}\in \mathrm{CR}^1_0(\hT)$ and $\widehat \beta_C \in S^1(\hT)/\mathbb R$ 
on the fine level, all   bounded by the left-hand side  $\lesssim 1$. 
The same argument shows   
\begin{equation}\label{marinipur}
	\hpRT^*= \nabla_{NC} \huCR^* - 1/2\, \left(\Pi_0 f\right) (\bullet -\midhT)
\end{equation}
for some $\huCR^*\in CR^1_0(\hT)$.
The remaining term $a(\hpRT^*,\qRT-\hqRT)$  equals 
\[
\int_\Omega \hpRT^*\cdot \nabla_{NC} \widehat{v_{CR}} \, dx+
\frac12 \int_\Omega\hpRT^*  \cdot (x -\midhT) (\Pi_0-1)\ddiv \hqRT   dx.
\]
This, the representation \eqref{aforementionedrepresentation} of $ \qRT-\hqRT$, and an
integration by parts show 
\begin{align*}
& \int_\Omega \hpRT^*\cdot \nabla_{NC} \widehat{v_{CR}} \, dx
	  = \int_\Omega \nabla_{NC} \huCR^*\cdot  \nabla_{NC} \hvCR  dx \\
	&  = \int_\Omega \nabla_{NC} \huCR^*\cdot( \qRT-\hqRT ) dx 
	 =\int_{\Omega'}   \huCR^*  \ddiv ( \hqRT-\qRT ) dx.
\end{align*}
Therein, $\Omega'$ is the interior of the $\bigcup (\T\setminus\hT)$, the union of the elements in $\T\setminus\hT$.
Since $ \ddiv ( \hqRT-\qRT )=(1-\Pi_0)\ddiv \hqRT$ is $L^2$ perpendicular to $P_0(\T)$ (and so
vanishes on $\T\cap\hT$ outside of $\Omega'$), 
a discrete Poincare inequality proves
that this is bounded from above by  
$\lesssim ||  h_\T   \nabla_{NC}   \huCR^*  ||_{L^2(\Omega')} $.
% Since $\hpRT^*  \cdot (x -\midhT)=: -(\Pi_0 f/2) \, s^2(\hT)$ for $s^2(\hT) \approx h^2_{\hT}\le h^2_\T$, 
Since  $(1-\Pi_0)\ddiv \hqRT$ vanishes outside of $\Omega'$ and has a
bounded $L^2$ norm,  the second integral reads
\[
\frac12 \int_\Omega\hpRT^*  \cdot (x -\midhT) (\Pi_0-1)\ddiv \hqRT   dx
\lesssim \|  h_\T \hpRT^* \|_{L^2(\Omega')}.
\]
The combination of the three previously displayed formulas and a triangle inequality lead to
\begin{align*}
a(\hpRT^*,\qRT-\hqRT)&\lesssim  \|  h_\T  \hpRT^* \|_{L^2(\Omega')}
+ ||  h_\T   \nabla_{NC}   \huCR^*  ||_{L^2(\Omega')} \\
&\lesssim  \|  h_\T  \hpRT^* \|_{L^2(\Omega')}
+ ||  h_\T  (  \hpRT^*- \nabla_{NC}   \huCR^* )  ||_{L^2(\Omega')}.
\end{align*}
The representation \eqref{marinipur} shows that the last term is equal to
\[
	1/2 \,  ||  h_\T  \, (\Pi_0 f)\, (\bullet -\midhT)  ||_{L^2(\Omega')}
\lesssim  ||  h_\T^2 \,\ddiv p_{RT}  ||_{L^2(\Omega')}.
\]
An inverse estimate for  $ \pRT$ on any $K\in\T\setminus\hT$ leads to 
\[
 \|  h_\T^2  \ddiv p_{RT} \|_{L^2(\Omega')} \lesssim  \|  h_\T  p_{RT} \|_{L^2(\Omega')}.
\]
A triangle inequality plus  $||  h_\T ||_{L^\infty(\Omega')}\lesssim 1$ prove
\[
 \|  h_\T  \hpRT^* \|_{L^2(\Omega')}\lesssim  
  \|  h_\T  \pRT \|_{L^2(\Omega')}
 +||  \hpRT^*-p_{RT}  ||_{L^2(\Omega')}.
\]
The combination of the above estimates (i.e.\ \eqref{eq:est:LHS},\eqref{eq:biA:CSI} and the three previously displayed formulas) shows that
\begin{eqnarray}\label{referenceestimate}
&&\| \hpRT-\pRT\|_{H(\ddiv,\Omega)} +\| \huRT-u_{RT} \|_{L^2(\Omega)}\\ 
\nonumber
&& \lesssim   ||  \hpRT^*-p_{RT}  ||_{L^2(\Omega)}+ \|  h_\T  \pRT \|_{L^2(\Omega')}
+ \|\widehat{\Pi_0} f-\Pi_0 f\|_{L^2(\Omega)}.
\end{eqnarray}
The $L^2$ orthogonal projection $\Pi_0$ (resp. $\widehat{\Pi_0}$) with respect to $\T\in\TO$
 (resp. its refinement $\hT\in\TO(\T)$) leads to the data approximation term  
 \[
 \|\widehat{\Pi_0} f-\Pi_0 f\|_{L^2(\Omega)}^2= \mu^2(\T)-\mu^2(\hT).
 \]
The combination of this with \eqref{stabilityintermediate} and \eqref{referenceestimate} proves (A3) in the sharper form
\[
\delta^2(\T,\hT)\le \| \hpRT-\pRT\|_{H(\ddiv,\Omega)}^2 +\| \huRT-u_{RT} \|_{L^2(\Omega)}^2
\lesssim \eta^2(\T,\T\setminus\hT)+  \mu^2(\T)-\mu^2(\hT).
\]

The proof of (A4) recalls the  $L^2$ quasiorthogonality of the flux errors of
\cite[Thm 3.2]{LCMHJX} or 
\cite[Lemma 4.3 and (4.4)]{CR09} in the form
\[
 \| p_{\ell+1} -p_{\ell}\|_{L^2(\Omega)}^2
 +\| p -p_{\ell+1}\|_{L^2(\Omega)}^2
 - \| p -p_{\ell}\|_{L^2(\Omega)}^2 
  \lesssim 
  \| p -p_{\ell+1}\|_{L^2(\Omega)}\, \osc(f_{\ell+1},\T_\ell).
\]
The mixed FEM fixes the divergence of the flux approximations, $-\ddiv p_\ell=\Pi_\ell f=: f_\ell$, and their orthogonality
\[
 \| f_{\ell+1} -f_{\ell}\|_{L^2(\Omega)}^2
 +\| f -f_{\ell+1}\|_{L^2(\Omega)}^2
 - \| f -f_{\ell}\|_{L^2(\Omega)}^2 =0
\]
leads (for all $\ell\in\mathbb{N}$) in the aforementioned  $L^2$ quasiorthogonality to 
\[
 \| p_{\ell+1} -p_{\ell}\|_{H(\ddiv,\Omega)}^2
 +\| p -p_{\ell+1}\|_{H(\ddiv,\Omega)}^2
 - \| p -p_{\ell}\|_{H(\ddiv,\Omega)}^2 
  \lesssim 
  \| p -p_{\ell+1}\|_{L^2(\Omega)}\, \osc(f_{\ell+1},\T_\ell).
\]
For any $0<\varepsilon$ with $\varepsilon \LdRel<1$ and the multiplicative constant $C\approx 1$ hidden in the notation $\lesssim$ 
the sum of those estimates results for any $\ell,m\in\mathbb{N}_0$ in 
\begin{align}\label{proofofa4inmfem}
	\begin{aligned}
		\sum_{k=\ell}^{\ell+m} \| p_{k+1} -p_{k}\|_{H(\ddiv,\Omega)}^2
		& \le  \| p -p_{\ell}\|_{H(\ddiv,\Omega)}^2 +   
\varepsilon/\LdRel \sum_{k=\ell}^{\ell+m-1}  \| p -p_{k+1}\|_{L^2(\Omega)}^2 \\
&\quad 
+ C^2\LdRel/\varepsilon \sum_{k=\ell}^{\ell+m}  \osc^2(f_{k+1},\T_k). 
\end{aligned}
\end{align}
For a sequence of uniformly refined meshes $\hT$, the discrete reliability (A3) leads  
to the reliability  of \cite{ccdpas2015}, 
\[
	\| p -p_{\ell}\|_{H(\ddiv,\Omega)}^2 \le \LdRel\,\eta_\ell^2:=\eta^2(\T_\ell)
 \quad\text{for all }\ell\in \mathbb{N}_0.
\]
The oscillation  $\osc(f_{k+1},\T_k)=\| h_\ell(f_{k+1}-f_k)\|_{L^2(\Omega)}$ is bounded by 
$\| h_\ell \|_{L^\infty(\Omega)}\, \| f_{k+1}-f_k\|_{L^2(\Omega)}$. 
With $h_{\max} :=  \| h_0 \|_{L^\infty(\Omega)} \lesssim 1$,
the $L^2$ orthogonality of the integrants shows
\[
 \sum_{k=\ell}^{\ell+m}  \osc^2(f_{k+1},\T_k)\le h_{\max} \, \| f_{\ell+m+1}-f_\ell\|_{L^2(\Omega)}^2
 \le h_{\max} \, \| f -f_\ell\|_{L^2(\Omega)}^2.
\]
The combination of the previous estimates with \eqref{proofofa4inmfem} leads to the quasiorthogonality
(A4) in the form
\[
\sum_{k=\ell}^{\ell+m} \delta_{k,k+1}^2
\le \LdRel\eta_\ell^2+ \varepsilon \, \sum_{k=\ell+1}^{\ell+m}\eta_k^2
+C \LdRel h_{\max}/\varepsilon \, \mu^2(\T_\ell).
\]
This is ($A4_\varepsilon$) with $\LqoE:=\max\{ \LdRel, C \LdRel h_{\max}/\varepsilon\}$ 
for any $\varepsilon>0$. This and (A1)-(A2) imply (A4) owing to \cref{thm:A4eA4}.
The remaining details are omitted. 
\qquad \endproof

\section{Application to least-squares FEM}
\label{s:appllsfem}
\newcommand{\uls}{u_{\operatorname{LS}}}
\newcommand{\pls}{p_{\operatorname{LS}}}
\newcommand{\qls}{q_{\operatorname*{LS}}}
\newcommand{\vls}{v_{\operatorname*{LS}}}
\newcommand{\sigmals}{\sigma_{\operatorname{LS}}}
\newcommand{\tauls}{\tau_{\operatorname*{LS}}}
The div least-squares formulation  \cite{BG09} of the Poisson model example of the previous section seeks 
the minimizer  $(p,u)$ of the functional
\begin{equation*}
	\operatorname{LS}(f; q,v):=
	\|f+\operatorname{div}q\|_{L^2(\Omega)}^2+\|q -\nabla v\|_{L^2(\Omega)}^2
\end{equation*}
amongst $(q,v)\in H(\ddiv,\Omega)\times H^1_0(\Omega)$. The functional $\operatorname{LS}(f;\bullet)$ is indeed a natural a~posteriori error estimator. Given any admissible triangulation $\T\in\TO$,
the least-squares FEM seeks the minimizer
$(\pls,\uls)$ of $\operatorname{LS}(f;\bullet)$ in the discrete subspace $RT_0(\T)\times S^1_0(\T)$. 
This leads in \cite{CCP-lsfem} 
to the alternative a~posteriori error estimate with
\begin{eqnarray}  \label{eq:estimatorls}
\tilde \eta^2(\T,K)&:=& \|(1-\Pi_0)\pls \|^2_{L^2(K)}
+|K|^{{1}/{2}} \hspace{-2mm}
\sum_{E\in\E(K)} \hspace{-1mm} \| [ \pls  ]_E \cdot \tau_E \|_{L^2(E)}^2\\  \nonumber
&& + |K|^{{1}/{2}}\hspace{-2mm} \sum_{E\in\E(K)\setminus\E(\partial\Omega)}\hspace{-1mm}
\| [ \partial \uls /\partial \nu_E]_E \|_{L^2(E)}^2, \\
\mu^2(K)&:=&|| f-\Pi_0 f||^2_{L^2(K)}\quad\text{for any }K\in\T.  \label{eq:muls}
\end{eqnarray}
Given a refined triangulation $\hT\in\TO(\T)$ with discrete solutions 
$(\widehat{\pls},\widehat{\uls})$, the distance 
\[
\delta^2(\T,\hT):= \operatorname{LS}(f; \pls,\uls)- \operatorname{LS}(f; \widehat{\pls},\widehat{\uls})
=\operatorname{LS}(0; \widehat{\pls}-\pls,\widehat{\uls}-\uls)\]
is equivalent to the norm of the difference $(\widehat{\pls}-\pls, \widehat{\uls} -\uls)$ of the two discrete solutions
in $H(\ddiv,\Omega)\times H^1_0(\Omega)$ \cite{BG09}.

\begin{theorem}[A1--A4] \label{thm:lsfem}
The estimators and distance function satisfy (A1)-(A4) and (B2) for 
$\mathcal{R}(\T,\hT):=\T\setminus\hT$, $\Lref=1=\LsubAdd$, and (QM). 
\end{theorem}

Since the estimator is reliable and efficient, for the discrete solution $(\pls,\uls)$ with respect to 
$\T\in\TO$, 
\[
\sigma(\T)\approx \| p-\pls\|_{H(\ddiv,\Omega)}+\| u-\uls\|_{H^1(\Omega)},
\] 
the optimal rates of the estimators is equivalent to the optimal rates of the errors in terms of nonlinear approximation 
classes with respect to the natural norms in $H(\ddiv)\times H^1$ of the least-squares FEM.

\textit{Proof of \cref{thm:lsfem}.}
The proofs are essentially contained in  \cite{CCP-lsfem}. The axioms (A1)-(A2) are standard and (A3) follows 
from
\[
\operatorname{LS}(f; \pls,\uls)\lesssim \eta^2(\T,\T\setminus\hT)+\mu^2(\T) 
+ \operatorname{LS}(f; \widehat{\pls},\widehat{\uls})
\]
(this is \cite[p 59, line 24]{CCP-lsfem} in different notation) and another reliability estimate 
$ \operatorname{LS}(f; \widehat{\pls},\widehat{\uls})\approx \sigma^2(\hT)$ from \cite[Thm 3.1]{CCP-lsfem}.
Notice that $\LhdRel$ does {\em not} need to be small (at least for coarse meshes according to the Remark \ref{rem:appLS}) and hence \cref{thm:qmono} cannot be applied to ensure (QM) in general. On the other hand, any conforming discretization reduces the least-squares functions and so 
\[
\sigma^2(\hT)\approx  \operatorname{LS}(f; \widehat{\pls},\widehat{\uls})\le 
 \operatorname{LS}(f; {\pls},{\uls})\approx \sigma^2(\T)
\]
immediately leads to (QM).
The same argument plus the reliability of \cite[Theorem 3.1]{CCP-lsfem} prove (A4) even in the sharper form of an orthogonality. The remaining details are omitted.
\qquad \endproof

\begin{remark}\label{rem:appLS}
A detailed analysis of  \cite{CCP-lsfem} (beyond this paper) with reduced elliptic regularity
suggests that $\LhdRel \le C(\epsilon) \, h_{\max}^{1/2+\epsilon} $ for small $\epsilon>0$ (depending on the interior angles of the domain) and some constant $C(\epsilon)$. Hence  $\LhdRel $ tends to zero as the maximal
mesh-size $h_{\max}$ tends to zero and so \cref{thm:qmono} is applicable for sufficiently fine meshes.
\end{remark}
\begin{remark}
	The analysis also allows optimal convergence rates for modified estimators such as
	\[
		\eta^2(K):= \abs{K} \norm{\D p_{LS}}{L^2(K)} + \abs{K}^{1/2} \norm{[p_{LS}-\nabla u_{LS}]_{\partial K}}{L^2(\partial K)}
	\]
	with $[p_{LS}-\nabla u_{LS}]_{\partial K}:= \restrict{(p_{LS}-\nabla u_{LS})}{K}$ along $E \in \E(\partial \Omega)$ with $K= \bar \omega_E$. This estimator is close to the least-squares functional estimators, but not equivalent.%, cf.\ \cite{lsfemBC} for details of the latter related algorithm. 
\end{remark}

%% file: axioms_safem.bbl
\begin{thebibliography}{CKNS08}

\bibitem[AB85]{ArnoldBrezzi}
D.~N. Arnold and F.\ Brezzi.
\newblock Mixed and nonconforming finite element methods: implementation,
  postprocessing and error estimates.
\newblock {\em RAIRO, Mathematical Modelling and Numerical Analysis},
  19(1):7--32, 1985.

\bibitem[Alo96]{Alonso96}
A.\ Alonso.
\newblock Error estimators for a mixed method.
\newblock {\em Numer. Math.}, 74:385--395, 1996.

\bibitem[BC]{lsfemBC}
P.\ Bringmann and C.\ Carstensen.
\newblock An adaptive least-squares {FEM} for the {S}tokes equations with
  optimal convergence rates.
\newblock Submitted to Numer.\ Math., 2015.

\bibitem[BCS]{BCStarke}
P.\ Bringmann, C.\ Carstensen, and G.\ Starke.
\newblock An adaptive least-squares {FEM} for linear elasticity with optimal
  convergence rates.
\newblock In preparation.

\bibitem[BDdV04]{BDD04}
P.\ Binev, W.\ Dahmen, and R.~de~Vore.
\newblock Adaptive finite element methods with convergence rates.
\newblock {\em Numer. Math.}, 97:219--268, 2004.

\bibitem[BdV04]{BD04}
P.\ Binev and R.~de~Vore.
\newblock Fast computation in adaptive tree approximation.
\newblock {\em Numer. Math.}, 97:193--217, 2004.

\bibitem[BG09]{BG09}
Pavel~B.\ Bochev and Max~D. Gunzburger.
\newblock {\em Least-Squares Finite Element Methods}.
\newblock Springer, 2009.

\bibitem[BM08]{mfemBeckerMao08}
R.\ Becker and S.\ Mao.
\newblock An optimally convergent adaptive mixed finite element method.
\newblock {\em Numer. Math.}, 111:35--54, 2008.

\bibitem[Car97]{CC-MC-97}
C.\ Carstensen.
\newblock A posteriori error estimate for the mixed finite element method.
\newblock {\em Math.\ Comp.}, 66:465--476, 1997.

\bibitem[CFPP14]{CFP14}
C.~Carstensen, M.~Feischl, M.~Page, and D.~Praetorius.
\newblock Axioms of adaptivity.
\newblock {\em Comput.\ Methods Appl.\ Math.}, 67(6):1195--1253, 2014.

\bibitem[CHX09]{LCMHJX}
L.\ Chen, M.\ Holst, and J.\ Xu.
\newblock Convergence and optimality of adaptive mixed finite element methods.
\newblock {\em Math.\ Comp.}, 78(265):35--53, 2009.

\bibitem[CKNS08]{CKNS07}
J.~M.\ Cascon, C.\ Kreuzer, R.~H.\ Nochetto, and K.~G. Siebert.
\newblock Quasi-optimal convergence rate for an adaptive finite element method.
\newblock {\em SIAM J. Numer. Anal.}, 46(5):2524--2550, 2008.

\bibitem[CP15]{CCP-lsfem}
C.\ Carstensen and E.-J. Park.
\newblock Convergence and optimality of adaptive least squares finite element
  methods.
\newblock {\em SIAM J. Numer. Anal.}, 53:43--62, 2015.

\bibitem[CPS15]{ccdpas2015}
C.\ Carstensen, D.\ Peterseim, and A.~Schr\"oder.
\newblock The norm of a discretized gradient in ${H}(\textup{div})^*$ for a
  posteriori finite element error analysis.
\newblock {\em Numer. Math.}, 2015.

\bibitem[CR11]{CR09}
C.~Carstensen and H.~Rabus.
\newblock An optimal adaptive mixed finite element method.
\newblock {\em Math. Comp.}, 80(274):649--667, 2011.

\bibitem[GSS14]{GSStevcmam}
Dietmar Gallistl, Mira Schedensack, and Rob Stevenson.
\newblock A remark on newest vertex bisection in any space dimension.
\newblock {\em Computational methods in applied mathematics}, 14(3), 2014.

\bibitem[Mar85]{Marini}
L.~D. Marini.
\newblock An inexpensive method for the evaluation of the solution of the
  lowest order {R}aviart-{T}homas mixed method.
\newblock {\em SIAM J.\ Numer.\ Anal.}, 22:493--496, 1985.

\bibitem[Rab15]{safem2015}
H.\ Rabus.
\newblock Quasi-optimal convergence of {AFEM} based on separate marking --
  {Part I and II}.
\newblock {\em Journal of Numerical Mathematics}, 23(2):137--156, 157--174,
  2015.

\bibitem[Ste07]{Stev07}
R.~Stevenson.
\newblock Optimality of a standard adaptive finite element method.
\newblock {\em Foundations of Computational Mathematics}, 7(2):245--269, 2007.

\bibitem[Ste08]{Stev08}
R.~Stevenson.
\newblock The completion of locally refined simplicial partitions created by
  bisection.
\newblock {\em Math. Comp.}, 77:227--241, 2008.

\end{thebibliography}
